\DeclareFontFamily{U}{wncy}{}
\DeclareFontShape{U}{wncy}{m}{n}{%
   <5>wncyr5%
   <6>wncyr6%
   <7>wncyr7%
   <8>wncyr8%
   <9>wncyr9%
   <10>wncyr10%
   <11>wncyr10%
   <12>wncyr6%
   <14>wncyr7%
   <17>wncyr8%
   <20>wncyr10%
   <25>wncyr10}{}
\DeclareMathAlphabet{\cyrille}{U}{wncy}{m}{n}
\def\sh{
\setlength{\unitlength}{.5 pt}
\begin{picture}(40,20)
\put(10,2){\line(1,0){20}}
\put(10,2){\line(0,1){10}}
\put(20,2){\line(0,1){10}}
\put(30,2){\line(0,1){10}}
\end{picture}}

\newcommand{\dlangle}{\langle\!\langle}
\newcommand{\drangle}{\rangle\!\rangle}

\documentclass[a4paper,11pt]{article}

\usepackage{latexsym}
\usepackage{amsmath,amsfonts,amssymb,amsthm}
\newtheorem{thm}{Theorem}[section]

\newtheorem{cor}[thm]{Corollary}
\newtheorem{prop}[thm]{Proposition}
\newtheorem{conj}[thm]{Conjecture}
\newtheorem{rem}[thm]{Remark}

\input epsf

\title{Arithmetic properties related to the shuffle-product}
\author{Roland Bacher
}
\date{}
\begin{document}
\maketitle

{\sl Abstract\footnote{Keywords: Shuffle product, formal power
  series, rational fraction, algebraic power series,
quadratic form, automaton sequence, Math. class: 11B85, 11E08, 11E76}: 
Properties of the shuffle product
in positive characteristic suggest to consider a $p-$homogeneous
form $\sigma:\overline{\mathbb F_p}\dlangle X_1,\dots,X_k\drangle
\longrightarrow
\overline{\mathbb F_p}\dlangle X_1,\dots,X_k\drangle$ on the vector space
$\overline{\mathbb F_p}\dlangle X_1,\dots,X_k\drangle$ of formal power series in
$k$ free non-commuting variables. The form $\sigma$ preserves rational
elements in
$\overline{\mathbb F_p}\dlangle X_1,\dots,X_k\drangle$, 
algebraic series of 
$\overline{\mathbb F_p}[[X]]=\overline{\mathbb F}\dlangle X\drangle$ 
and induces a bijection 
on the affine subspace $1+\mathfrak m$ of formal power series with
constant coefficient $1$. Conjecturally, this bijection restricts to a
bijection of rational elements in $1+\mathfrak m\subset
\overline{\mathbb F_p}\dlangle X_1,\dots,X_k\drangle$, 
respectively algebraic
elements in $1+X\overline{\mathbb F_p}[[X]]$.}

\section{Introduction}

The aim of this paper is to present some properties and conjectures 
related to shuffle-products of power series in non-commuting 
variables. The shuffle product
$$A\sh B=\sum_{0\leq i,j} {i+j\choose i}\alpha_i\beta_j X^{i+j}$$
of two power series $A=\sum_{n=0}^\infty \alpha_nX^n,
B=\sum_{n=0}^\infty \beta_nX^n\in \mathbb K[[X]]$ in one variable over a 
commutative field $\mathbb K$ 
turns the set $\mathbb K^*+X\mathbb K[[X]]$
into a commutative group which is not isomorphic to the commutative 
group on $\mathbb K^*+X\mathbb K[[X]]$ associated to the
ordinary product of (multiplicatively) invertible
formal power series if $\mathbb K$ is of positive characteristic. 
Shuffle products of
rational (respectively algebraic) power series are rational
(respectively algebraic).
The shuffle product turns the affine subspace
$1+X\mathbb K[[X]]$ into a group which is isomorphic to an
infinite-dimensional $\mathbb F_p-$vector space if $\mathbb K$ is a
field of positive characteristic $p$.
Rational (respectively algebraic) elements in 
$1+X\mathbb K[[X]]$ (or more generally in $\mathbb K^*+X\mathbb
K[[X]]$) form thus a group with respect to the shuffle product
if $\mathbb K$ is of positive characteristic. 

The first interesting case is given by a subfield $\mathbb K\subset
\overline{\mathbb F_2}$ contained in the algebraic closure
of the field $\mathbb F_2$ with two elements.
The structure of the $\mathbb F_2-$vector
space induced by the shuffle product on $1+X\overline{\mathbb F_2}[[X]]$
suggests to consider the quadratic form 
$$\begin{array}{ll}\displaystyle
\sigma\left(\sum_{n=0}^\infty \alpha_nX^n\right)&
\displaystyle =
\sum_{n=0}^\infty \alpha_{2^n}^2X^{2^{n+1}}+\sum_{0\leq i\leq j}
{i+j\choose i}\alpha_i\alpha_jX^{i+j}\\
&\displaystyle=
\alpha_0^2+\sum_{n=0}^\infty \alpha_{2^n}^2X^{2^{n+1}}+\sum_{0\leq i<j}
{i+j\choose i}\alpha_i\alpha_jX^{i+j}\ .\end{array}$$
The quadratic form $\sigma:\overline{\mathbb F_2}[[X]]\longrightarrow
\overline{\mathbb F_2}[[X]]$ thus defined preserves the vector
space of rational or algebraic power series. 
It induces a bijection of infinite order 
on the affine subspace $1+X\overline{\mathbb F_2}[[X]]$.
Orbits are either infinite or of cardinality a power of two.
Conjecturally, the inverse bijection $\sigma^{-1}$ of the set 
$1+X\overline{\mathbb F_2}[[X]]$ preserves also rational elements
and algebraic elements. We present experimental evidence
for this conjecture.
An analogous construction yields a homogeneous $p-$form 
(still denoted)
$\sigma:\overline{\mathbb F_p}[[X]]\longrightarrow 
\overline{\mathbb F_p}[[X]]$ with similar properties for $p$ 
an arbitrary prime. 

In a second part of the paper, starting with Section 
\ref{sectionpows},we recall the definition of the 
shuffle product for elements in the vector space
$\mathbb K\dlangle X_1,\dots,X_k\drangle$
of formal power series in free 
non-commuting variables. The shuffle product preserves again
rational formal power series, characterised for instance 
by a Theorem of Sch\"utzenberger.
The $p-$homogeneous form $\sigma$ considered above has a natural
extension $\sigma:\overline{\mathbb
F_p}\dlangle X_1,\dots,X_k\drangle\longrightarrow \overline{\mathbb
F_p}\dlangle X_1,\dots,X_k\drangle$. 
This extension of $\sigma$ still preserves
rational elements and induces a bijection on $1+\mathfrak m$
where $\mathfrak m\subset \overline{\mathbb F_p}\dlangle X_1,\dots,X_k\drangle$ 
denotes the maximal ideal of formal power series without constant 
coefficient in $\overline{\mathbb F_p}\dlangle X_1,\dots,X_k\drangle$. 
Conjecturally, the map $\sigma$ restricts again to a bijection 
of the subset of rational elements in $1+\mathfrak m$.

\section{Power series in one variable}

We denote by $\mathbb K[[X]]$ the commutative algebra of formal power
series over a commutative field $\mathbb K$ with product
$$\left(\sum_{n=0}^\infty \alpha_nX^n\right)\left(\sum_{n=0}^\infty
\beta_nX^n\right)=\sum_{n,m=0}^\infty \alpha_n\beta_mX^{n+m}$$
given by the usual (Cauchy-)product 
extending the product of the polynomial subalgebra
$\mathbb K[X]\subset \mathbb K[[X]]$.
Its unit group $\mathbb K^*+X\mathbb K[[X]]$ consists of all
(multiplicatively) invertible series and decomposes as a direct 
product $\mathbb K^*\times(1+\mathfrak m)$
with $\mathfrak m=X\mathbb K[[X]]$ denoting the maximal ideal of 
the algebra $\mathbb K[[X]]$.

A subalgebra containing the field of constants 
$\mathbb K$ of $\mathbb K[[X]]$ 
is {\it rationally closed} if
it intersects the unit group $\mathbb K^*+X\mathbb K[[X]]$ in a subgroup. 
The {\it rational closure} of a subset $\mathcal S\subset \mathbb
K[[X]]$ is the smallest rationally closed subalgebra of 
$\mathbb K[[X]]$ which contains $\mathcal S$ and the ground-field
$\mathbb K$. 

The rational closure of $X$, called the {\it algebra of rational
  fractions in $X$} or the {\it rational subalgebra} of $\mathbb K[[X]]$, 
contains the polynomial subalgebra 
$\mathbb K[X]$ and is formed by all rational fractions of the form
$\frac{f}{g}$ with $f,g\in\mathbb K[X],
g\not\in\mathfrak m$. The expression $\frac{f}{g}$ of such a 
rational fraction is unique if we require $g\in 1+\mathfrak m$.

An element $y\in\mathbb K[[X]]$ is {\it algebraic} if it 
satisfies a polynomial identity
$P(X,y)=0$ for some polynomial $P\in\mathbb K[X,y]$. 
Algebraic series in $\mathbb K[[X]]$
form a rationally closed subalgebra containing all rational 
fractions.

\section{The shuffle product}\label{sectshufflpr}

The {\it shuffle product}, defined as
$$A\sh B=\sum_{n,m=0}^\infty {n+m\choose n}\alpha_n\beta_m X^{n+m}$$
for $A=\sum_{n=0}^\infty \alpha_n X^n,B=\sum_{n=0}^\infty \beta_n X^n
\in \mathbb K[[x]]$, yields an associative and commutative
bilinear product on the vector space $\mathbb K[[x]]$ of 
formal power series.
We call the corresponding algebra $(\mathbb K[[x]],\sh)$ the
{\it shuffle-algebra}. The {\it shuffle-group} is the 
associated unit group. Its elements are given by the set $\mathbb K^*
+X\mathbb K[[x]]$ underlying the multiplicative unit group and it
decomposes as a direct product $\mathbb K^*\times(1+X\mathbb K[[X]])$.

\begin{rem} \label{remisomalg}
Over a field $\mathbb K$ of characteristic zero, the map
$$\mathbb K[[X]]\ni\sum_{n=0}^\infty \alpha_nX^n
\longmapsto \sum_{n=0}^\infty 
n!\alpha_nX^n\in(\mathbb K[[X]],\sh)$$
defines an isomorphism of algebras between the usual (multiplicative)
algebra of formal power series and the shuffle algebra 
$(\mathbb K[[X]],\sh)$.
The shuffle product of ordinary generating series $\sum \alpha_n X^n$
corresponds thus to the ordinary product of exponential generating series 
(also called divided power series or Hurwitz series, see eg. \cite{C})
$\sum\alpha_n \frac{X^n}{n!}$. 
This shows in particular the identity
$(1-X)\sh\left(\sum_{n=0}^\infty n!X^n\right)=1$. The shuffle inverse
of a rational fraction is thus generally transcendental
in characteristic $0$.
\end{rem}

\begin{rem}\label{remamoins1inverse} 
The inverse for the shuffle product of 
$1-a\in 1+X\mathbb K[[x]]$ is given by
$$\sum_{n=0}^\infty a^{\sh^n}=1+a+a\sh a+a\sh a\sh
a+\dots$$
where $a^{\sh^0}=1$ and $a^{\sh^{n+1}}=a\sh a^{\sh^n}$ for $n\geq 1$.

The shuffle inverse of $1-a\in A+X\mathbb K[[X]]$ can 
be computed by the recursive formulae
$B_0=1,C_0=a,B_{n+1}=B_n+B_n\sh C_n,C_{n+1}=C_n\sh C_n=
a^{\sh^{2^{n+1}}}$ with
$B_n=\sum_{k=0}^{2^n-1}a^{\sh ^k}$ converging (quadratically)
to the shuffle-inverse 
of $1-a$.
\end{rem}

\begin{prop} \label{propshgrFp} 
The shuffle-group $1+X\mathbb K[[X]]$ is isomorphic to an
infinite-dimensional $\mathbb F_p-$vector-space if 
$\mathbb K$ is a field of positive characteristic $p$.
\end{prop}

\begin{cor} The shuffle-group $1+X\mathbb K[[X]]$ is
not isomorphic to the multiplicative group structure on $1+X\mathbb
K[[X]]$ if $\mathbb K$ is of positive characteristic.
\end{cor}

{\bf Proof of Proposition \ref{propshgrFp}} We have
$$A^{\sh^p}=\sum_{0\leq i_1,i_2,\dots,i_p}{i_1+i_2+\dots+i_p\choose
  i_1,i_2,\dots,i_p}\alpha_{i_1}\alpha_{i_2}\cdots 
\alpha_{i_p}X^{i_1+\dots+i_p}  \ .$$
for $A=\sum_{n=0}^\infty \alpha_nX^n\in\mathbb K[[X]]$ where 
${i_1+i_2+\dots+i_p\choose
  i_1,i_2,\dots,i_p}=\frac{(i_1+\dots +i_p)!}{i_1!\cdots i_p!}$.
Two summands differing by a cyclic permutation 
of indices $(i_1,i_2,\dots,i_p)\longmapsto
(i_2,i_3,\dots,i_p,i_1)$ yield the same contribution to $A^{\sh^p}$. 
Over a field $\mathbb K$ of positive characterstic $p$ we 
can thus restrict the summation to $i_1=i_2=\dots=i_p$.
Since ${ip\choose i,i,\dots,i}=\frac{(ip)!}{(i!)^p}\equiv 0\pmod p$
except for $i=0$, we have $A^{\sh^p}=\alpha_0^p$ for
$A=\sum_{n=0}^\infty
\alpha_nX^n\in\mathbb K[[X]]$. This implies the result .
\hfill$\Box$

\begin{rem} Proposition \ref{propshgrFp} follows also easily 
from Satz 1 in \cite{H} where a different proof is given.
\end{rem}

\begin{prop} \label{proprat} Shuffle products of rational power series are
rational.
\end{prop}

{\bf Proof} Suppose first
$\mathbb K$ of characteristic zero. 
The result is obvious for the shuffle product
of two polynomials. Extending $\mathbb K$ to its algebraic closure,
decomposing into simple fractions 
and using bilinearity, it is enough to consider
shuffle products of the form
$X^h\sh\left(\sum_{n=0}^\infty n^k\alpha_nX^n\right)=
\sum_{n=0}^\infty {n+h\choose h}n^k\alpha^nX^{n+h}$ 
which are obviously rational
and shuffle products of the form
$$\left(\sum_{n=0}^\infty n^h\alpha^nX^n\right)\sh
 \left(\sum_{n=0}^\infty n^k\beta^nX^n\right)=
\sum_{0\leq m\leq n} {n\choose m}m^h(n-m)^k\alpha^m\beta^{n-m}X^n$$
which are evaluations at $y=\alpha,z=\beta$ of
$$\left(y\,\frac{\partial}{\partial y}\right)^h\left(z\,
\frac{\partial}{\partial z}\right)^k\left(\frac{1}{1-(y+z)X}\right)
$$
and are thus rational for $\mathbb K$ of characteristic zero.

In positive characteristic, one can either consider suitable lifts into 
integer rings of fields of characteristic zero or deduce it as a
special case of Corollay \ref{corratclosed}.\hfill$\Box$

\begin{rem} The proof of proposition \ref{proprat} implies easily 
analyticity of shuffle products of analytic power
series (defined as formal power series with strictly positive
convergence radii) if 
$\mathbb K\subset \mathbb C$ or $\mathbb K\subset \hat{\overline{
\mathbb Q_p}}$.
\end{rem}

\begin{prop} \label{propalg} Shuffle products of algebraic series in 
$\overline{\mathbb F_p}[[X]]$ are algebraic.
\end{prop}

{\bf Sketch of Proof} A Theorem
of Christol (see Theorem 12.2.5 in \cite{AS}) states that
the coefficients of an algebraic series over
$\subset \overline{\mathbb F_p}$ 
define a $q-$automatic sequence
with values in $\mathbb F_q$ for some power $q=p^e$ of $p$.
Given a formal power series $C=\sum_{n=0}^\infty \gamma_nX^n
\in \overline{\mathbb F_p}[[X]]$, we
denote by $C_{k,f}$ the formal power series $\sum_{n=0}^\infty 
\gamma_{k+nq^f}X^n$.

The result follows then from the observation that the series 
$(A\sh B)_{k,f}$ are linear combination of 
$A_{k_1,f}\sh B_{k_2,f}$ and span thus a finite-dimensional 
subspace of $\overline{\mathbb F_p}[[X]]$ for algebraic 
$A,B\in\overline{\mathbb F_p}[[X]]$.\hfill$\Box$ 

Propositions \ref{propshgrFp} and \ref{proprat} 
(respectively \ref{propshgrFp} and \ref{propalg})
imply immediately the following result:

\begin{cor} Rational (respectively algebraic)
  elements of the shuffle-group 
$1+X\mathbb K[[X]]$ form a subgroup for $\mathbb
K\subset\overline{\mathbb F_p}$.
\end{cor}

\begin{rem} \label{remshfflegrpnoalgelts}
A rational fraction 
$A\in 1+X\mathbb C[[X]]$
has a rational inverse for the shuffle-product if and only if
$A=\frac{1}{1-\lambda X}$ with
$\lambda\in\mathbb C$. (Idea of proof: Decompose two rational
series $A,B$ satisfying $A\sh B=1$ into simple fractions and
compute $A\sh B$
using the formulae given in the proof of Proposition \ref{proprat}.)
\end{rem}

\section{A quadratic form}

The identity $A\sh A=\alpha_0^2$
for $A=\sum_{n=0}^\infty \alpha_nX^n\in\overline{\mathbb F_2}[[X]]$
(see the proof of Proposition \ref{propshgrFp}) suggests to consider
the quadratic map 
$${\mathbb K}[[X]]\ni A=
\sum_{n=0}^\infty \alpha_nX^n\longmapsto \sigma(A)=\alpha_0^2+
\sum_{n=1}^\infty \beta_nX^n\in {\mathbb K}[[X]]\subset
\overline{\mathbb F_2}[[X]]$$
defined by 
$$\left(\sum_{n=0}^\infty \tilde\alpha_nX^n\right)\sh
\left(\sum_{n=0}^\infty \tilde\alpha_nX^n\right)=\tilde\alpha_0^2+
2\sum_{n=0}^\infty \tilde\beta_nX^n$$
where $\tilde\alpha_n$ and $\tilde \beta_n$ are lifts into 
suitable algebraic integers of $\alpha_n,\beta_n\in 
\mathbb K\subset\overline{\mathbb F_2}$.

For $A=\sum_{n=0}^\infty \alpha_nX^n$, we get
$$\sigma(A)=\alpha_0^2+\sum_{n=1}^\infty \frac{1}{2}{2n\choose n}
\alpha_n^2X^{2n}+\sum_{0\leq i<j}{i+j\choose
  i}\alpha_i\alpha_jX^{i+j}$$
and ${2n\choose n}\equiv 2\pmod 4$ if and only if $n$ is a power of
$2$.
This yields the formula
$$\sigma(A)=\alpha_0^2+\sum_{n=0}^\infty \alpha_{2^n}^2X^{2^{n+1}}+
 \sum_{0\leq i<j}{i+j\choose i}\alpha_i\alpha_jX^{i+j}\ .$$

\begin{prop} \label{propsigmaratalg}
The formal power series $\sigma(A)$ is rational
  (respectively algebraic) if $A\in\overline{\mathbb F_2}[[X]]$
is rational (respectively algebraic).
\end{prop}

The statement of this proposition in the case of a rational series 
is a particular case of Proposition \ref{propphomXk}.

Proposition \ref{propsigmaratalg} can be proven by modifying slightly
the arguments used in the proof  of Propositions \ref{proprat} and 
\ref{propalg} and by applying them to a suitable integral 
lift $\tilde A\in \overline{\mathbb Q}[[X]]$ of $A$. \hfill $\Box$

Finally, one has also the following result whose easy proof is left to
the reader:

\begin{prop} The quadratic form $A\longmapsto \sigma(A)$ commutes 
with the Frobenius
map $A\longmapsto A^2$.
\end{prop}

\subsection{The main conjecture}

\begin{prop} \label{propbij} The quadratic form $A\longmapsto \sigma(A)$
induces a bijection on the affine subspace $1+X\mathbb K[[X]]$ 
for a subfield $\mathbb K\subset\overline{\mathbb  F_2}$.
\end{prop}

\begin{rem} \label{remnotinjsurj}
Omitting the restriction to $1+X{\mathbb
    K}[[X]]$, the quadratic form $\sigma$ is neither 
surjective nor injective: One has 
$\sigma^{-1}(X)=\emptyset$ and $\sigma(A)=0$ if $A\in X^3\overline{
\mathbb K}[[X^2]]$. (The example for non-injectivity is related to
the easy observation that $\sigma(A)=0$ if and only if $\sigma(1+A)=1+A$
for $A\in\overline{\mathbb F_2}[[X]]$.)
\end{rem}

{\bf Proof of Proposition \ref{propbij}} This follows from the identity
$$\sigma(A)-\sigma(B)=(\alpha_n-\beta_n)X^n+X^{n+1}\overline{\mathbb
  F_2}[[X]]$$
if $A=1+\sum_{n=1}^\infty \alpha_nX^n,B=1+\sum_{n=1}^\infty \beta_nX^n$
coincide up to $X^{n-1}$ (ie. if $\alpha_j=\beta_j$ for
$j=1,\dots,n-1$).
\hfill$\Box$

Experimental evidence (see Sections \ref{sectratexples},
\ref{sectiterexples}
and \ref{sectalgcplproof} for a few exemples)
suggests the following conjecture:

\begin{conj} \label{conjX}
If $A\in 1+\overline{\mathbb F_2}[[X]]$ is rational
  (respectively algebraic) then its preimage $\sigma^{-1}(A)\in
1+\overline{\mathbb F_2}[[X]]$ is rational (respectively algebraic).
\end{conj}

This conjecture, in the case of rational power series, is a 
particular case of Conjecture \ref{conjratX_k} 
(which has, to my knowledge, no algebraic analogue).

\begin{rem}
There is perhaps some hope for proving 
this conjecture in the rational case 
using the formulae of the proof of Proposition  
\ref{proprat}: Considering integral lifts into suitable
algebraic integers and assuming a bound on the degrees
of the numerator and denominator of $\sigma^{-1}(A)$ (for 
rational $A\in 1+X\overline{\mathbb F_2}[[X]]$) one gets a system of
algebraic equations whose reduction modulo $2$ should have a 
solution.
\end{rem}

\subsection{Orbits in $1+X\overline{\mathbb F_2} [[X]]$ under $\sigma$}

The purpose of this Section is to describe a few properties of
the bijection defined by $\sigma$ on $1+X\overline{\mathbb F_2}[[X]]$.

\begin{prop} \label{propsigmaorbit}
(i) The orbit of $A\in1+X\overline{\mathbb F_2}[[X]]$ is infinite if it
involves a monomial of the form $X^{2^k}$.

\ \ (ii) The orbit of a polynomial 
$A\in 1+X\overline{\mathbb F_2}[X]$ is finite if it involves
no monomial of the form $X^{2^k}$.

\ \ (iii) 
The cardinal of every finite orbit in $1+X\overline{\mathbb F_2}[[X]]$
of $\sigma$ is a power of $2$. 
\end{prop}

\begin{rem} \label{remorbinf} (i) All elements of the form
$1+X^3\overline{\mathbb F_2}[[X^2]]$ are fixed by $\sigma$, cf. 
Remark \ref{remnotinjsurj}.

\ \ (ii) The algebraic function $A=1+\sum_{n=0}^\infty X^{3\cdot 4^n}$
(satisfying the equation $A+A^4+X^3=0$) contains no monomial
of the form $X^{2^k}$ and has infinite orbit under
$\sigma$. I ignore if the affine subspace
$1+X\overline{\mathbb F_2}[[X]]$ contains an infinite orbit
formed by rational fractions without monomials of the form $X^{2^k}$.
\end{rem}

{\bf Proof of Proposition \ref{propsigmaorbit}} Associate to 
$A=1+\sum_{n=1}^\infty \alpha_nX^n\in \overline{\mathbb F_2}[[X]]$ 
the auxiliary series $P_A=\sum_{n=0}^\infty
\alpha_{2^n}t^n\in\mathbb \overline{\mathbb F_2}[[t]]$. 
It is easy to check 
that $P_{\sigma^k(A)}=(1+t)^kP_A$ for all $k\in\mathbb Z$.
This implies assertion (i).

Consider a polynomial $A$ containing only coefficients of degree
$<2^n$ and no coefficient of degree a power of $2$. The formula
for $\sigma(A)$ shows that $\sigma(A)$ satisfies the same conditions.
This implies that the orbit of $A$ under $\sigma$ is finite and proves
assertion (ii).

If $A\in 1+\overline{\mathbb F_2}[[X]]$ is such that $\sigma^{2^k}(A)
\equiv A\pmod{X^{N-1}}$, then 
$\sigma^{2^k}(A+X^N)=\sigma^{2^k}(A)+X^N\pmod{X^{N+1}}$.
This implies easily the last assertion.\hfill$\Box$

\subsection{A variation}

The series $P_A=\sum_{n=0}^\infty \alpha_{2^n}t^n$ associated
to an algebraic power series $A=\sum_{n=0}^\infty \alpha_nX^n\in
\overline{\mathbb F_2}[[X]]$
as in the proof of proposition \ref{propsigmaorbit}
is always ultimately periodic and 
thus rational. This implies algebraicity of $\sum_{n=0}^\infty
\alpha_{2^n}X^{2^{n+1}}$ for algebraic $\sum_{n=0}^\infty \alpha_n
X^n\in\overline{\mathbb F_2}[[X]]$. The properties of the quadratic form 
$$A=\sum_{n=0}^\infty \alpha_nX^n\longmapsto \tilde
\sigma(A)=\sum_{0\leq i\leq j}{i+j\choose i}\alpha_i\alpha_jX^{i+j}$$
with respect to algebraic elements in $\overline{\mathbb F_2}[[X]]$ 
should thus be somewhat similar to the properties of $\sigma$.
It particular $\tilde \sigma$  preserves 
algebraic series and induces a bijection on
$1+X\overline{\mathbb F_2}[[X]]$ 
which is of infinite order. Orbits are either infinite or finite and
the cardinality of a finite orbit is 
a power of $2$. Conjecture \ref{conjX} (if true), 
together with Proposition \ref{propalg}, would imply that 
$\tilde \sigma^{-1}(A)$
is algebraic for algebraic $A\in 1+X\overline{\mathbb F_2}[[X]]$.
Remark however that $\tilde \sigma(A)$ is in general not rational
for rational $A\in 1+X\overline{\mathbb F_2}[[X]]$: An easy
computation shows indeed that $\tilde \sigma(\frac{1}{1+X})=
1+\sum_{n=0}^\infty X^{2^n}$ which satisfies the algebraic 
equation $y+y^2+X=0$ but is irrational since coefficients of rational 
power series over (the algebraic closure of) finite fields 
are ultimately periodic. 
On the other hand,
$\tilde\sigma^{-1}(\frac{1}{1+X})$
is the irrational 
algebraic series $y=1+X+X^2+X^4+X^7+\dots\in\mathbb F_2[[X]]$ 
satisfying the equation $$X+(1+X+X^2)y+(1+X^2+X^4)y^3=0\ .$$

The quadratic map $\tilde \sigma$
behaves however better than $\sigma$ with respect to polynomials:
One can show easily that it induces a bijection of order 
a power of $2$ (depending on $n$) on polynomials of degree $<2^n$
in $1+X\overline{\mathbb F_2}[[X]]$.

\begin{rem} The definition of the quadratic forms $\sigma$ and
$\tilde \sigma$ suggests
to consider the quadratic form 
$\psi(\sum_{n=0}^\infty \alpha_n X^n)=\sum_{i\leq
  j}\alpha_i\alpha_jX^{i+j}$ of $\overline{\mathbb F_2}[[X]]$.
Using the fact that rational elements of $\overline{\mathbb F_2}[[X]]$
have ultimately periodic coefficients, it is not hard to show
that $\psi$ preserves rationality. It is also easy to show that $\psi$
induces a bijection on $1+X\overline{\mathbb F_2}[[X]]$.
However, the preimage $\psi^{-1}(1+X)\in\mathbb F_2[[X]]$ 
is apparently neither rational nor algebraic.
\end{rem}

\subsection{Algorithmic aspects}

The {\it integral Thue-Morse} function 
$\mathop{tm}(\sum_{j=0} \epsilon_j 2^j)=\sum_j \epsilon_j$ 
is defined as the digit sum of a natural binary integer 
$n=\sum_{j=0} \epsilon_j
2^j\in \mathbb N$. Setting $\mathop{tm}(0)=0$,
it can then be computed recursively by
$\mathop{tm}(2n)=\mathop{tm}(n)$
and $\mathop{tm}(2n+1)=1+\mathop{tm}(n)$. Kummer's equality
${i+j\choose i}\equiv
2^{\mathop{tm}(i)+\mathop{tm}(j)-\mathop{tm}(i+j)}
\pmod 2$ (which follows also from a Theorem of Lucas,
see page 422 of \cite{AS}), allows a fast computation 
of binomial coefficients modulo $2$.
We have thus 
$$\begin{array}{rcl}
\displaystyle \sigma(A)&\displaystyle =&
\displaystyle \alpha_0^2+\sum_{n=0}^\infty \alpha_{2^n}^2X^{2^{n+1}}+
\sum_{0\leq i<j}{i+j\choose i}\alpha_i\alpha_jX^{i+j}\\
&\displaystyle =
&\displaystyle  \alpha_0^2+\sum_{n=0}^\infty \alpha_{2^n}^2X^{2^{n+1}}+
\sum_{0\leq i<j,\ \mathop{tm}(i+j)=\mathop{tm}(i)+\mathop{tm}(j)}
\alpha_i\alpha_jX^{i+j}\end{array}$$
for $A=\sum_{n=0}\alpha_n X^n\in\mathbb \overline{\mathbb F_2}[[x]]$.
The last formula is suitable for computations.

The preimage $\sigma^{-1}(A)$ of $A\in 1+X\overline{\mathbb
    F_2}[[X]]$ can be computed iteratively as the unique 
fixpoint in $\overline{\mathbb F_2}[[X]]$  of the map
$$Z\longmapsto Z+A-\sigma(Z)\ .$$
Starting with an arbitrary initial value $Z_0$ (eg. with $Z_0=A$), 
the sequence $Z_0,Z_1,\dots,Z_{n+1}=Z_n+A-\sigma(Z_n),\dots
\subset \overline{\mathbb F_2}[[X]]$ converges quadratically
(roughly doubling the number of correct coefficients at each 
iteration) with limit the attractive fixpoint $\sigma^{-1}(A)$.

\subsubsection{Checking identities in the rational case} 
Define the degree of a non-zero
rational fraction $A=\frac{f}{g}\in\overline{\mathbb F_2}[[X]]$ with
$f\in\overline{\mathbb F_2}[X],g\in 1+\overline{\mathbb
    F_2}[X]$ coprime,
by $\mathop{deg}(A)=\mathop{max}(\mathop{deg}(f),\mathop{deg}(g))$.
Proposition \ref{propphomXk} and Remark \ref{remdimratfrac}
imply the equality
$$\mathop{deg}(\sigma(A))\leq 1+{\mathop{deg}(A)+2\choose 2}\ .$$
This inequality can be used to prove identities 
of the form $\sigma(A)=B$ involving two rational fractions
$A,B\in \overline{\mathbb
  F_2}[X]$ by checking equality of the first
$2+{\mathop{deg}(A)+2\choose 2}+\mathop{deg}(B)$ 
coefficients of the series $\sigma(A)$ and $B$.

\subsubsection{Checking identities in the algebraic
  case}\label{sectidalgc} 
Given a power series 
$A=\sum_{n=0}^\infty \alpha_nX^n\in\overline{\mathbb F_2}[[X]]$, 
we consider the power series 
$A_{k,f}=\sum_{n=0}^\infty\alpha_{k+n\cdot 2^f}X^n$
for $k,f\in\mathbb N$ such that $0\leq k<2^f$. The 
vector space $\mathcal K(A)$ (called the $2-$kernel of $A$,
see \cite{AS})
spanned by all series $A_{k,f}$ is finite-dimensional if and only if $A$ 
is algebraic and one has the inequality
$$\mathop{dim}(\mathcal K(\sigma(A)))\leq 1+{1+\mathop{dim}(\mathcal
  K(A))\choose 2}\ .$$
This inequality, together with techniques of \cite{B}, 
reduces the proof of equalities $\sigma(A)=B$ involving
algebraic series $A,B\in\overline{\mathbb F_2}[[X]]$ to the
equality among finite series developpements of sufficiently
high order $N$ (depending on combinatorial properties) of $A$ and
$B$. The typical value for $N$ is of order 
$2^{2+{1+\mathop{dim}(\mathcal
  K(A))\choose 2}}$ and is thus unfortunately of no practical 
use in many cases.

\subsection{Examples involving  rational fractions in $1+\mathbb F_2[[X]]$}\label{sectratexples}

\subsubsection{A few preimages of polynomials}
$\sigma^{-1}(1+X)=\frac{1}{1+X}$, 
$\sigma^{-1}((1+X)^3)=1+X+X^3$, 
$\sigma^{-1}((1+X)^5)=(1+X)^2(1+X+X^2)(1+X^2+X^3)$,
$\sigma^{-1}((1+X)^7)=\frac{1+X^3+X^6}{(1+X)^7}$,
$\sigma^{-1}((1+X)^9)=(1+X)^6(1+X+X^9)$,
$\sigma^{-1}(1+X+X^2)=1+X$,
$\sigma^{-1}(1+X^2+X^3)=\frac{1+X^2+X^3}{(1+X)^4}$,
$\sigma^{-1}(1+X+X^3)=\frac{1+X+X^2}{(1+X)^4}$,
$\sigma^{-1}(1+X+X^4)=1+X+X^2+X^3$,
$\sigma^{-1}(1+X^3+X^4)=\frac{1+X+X^2}{(1+X)^3}$,
$\sigma^{-1}(1+X+X^2+X^3+X^4)=\frac{1+X+X^3}{(1+X)^4}$,
$\sigma^{-1}(1+X+X^2+X^3+X^5)=(1+X)(1+X^3+X^4)$, 
$\sigma^{-1}(1+X+X^3+X^4+X^5)=(1+X+X^2+X^5+X^7)$,
$\sigma^{-1}(1+X^2+X^3+X^4+X^5)=(1+X+X^3)(1+X+X^4)$,
$\sigma^{-1}(1+X^2+X^5)=\frac{(1+X+X^2)(1+X+X^3)}{(1+X)^6}$,
$\sigma^{-1}(1+X+X^2+X^4+X^5)=\frac{(1+X+X^4)}{(1+X)^8}$,
$\sigma^{-1}((1+X+X^2)^3)=\frac{1+X^2+X^3}{(1+X)^7}$,
$\sigma^{-1}((1+X)(1+X+X^2)=(1+X)(1+X+X^2)$, 
$\sigma^{-1}((1+X)^2(1+X+X^2))=(1+X+X^2)$, 
$\sigma^{-1}((1+X)^3(1+X+X^2)=\frac{1+X^3+X^4}{(1+X)^6}$,
$\sigma^{-1}((1+X)^4(1+X+X^2)=\frac{1+X+X^4+X^6+X^7}{(1+X)^8}$,

These examples suggest the following conjecture: 

\begin{conj} For $P\in 1+X\mathbb F_2[X]$ a polynomial of degree $\leq
  2^k$, we have $\sigma^{-1}(P)=\frac{Q_P}{(1+X)^{\alpha_P}}$
with $0\leq \alpha_P\leq 2^k$ and $Q_P\in 1+X\mathbb F_2[X]$
a polynomial of degree $<2^k$.
\end{conj} 

\subsubsection{A few examples of rational fractions}
$\sigma^{-1}\left(\frac{1}{(1+X)^3}\right)=\frac{(1+X)^2(1+X+X^4)}
{(1+X+X^2)^4}$, 
$\sigma^{-1}\left(\frac{1}{1+X+X^2}\right)=\frac{(1+X)^3}
{1+X^3+X^4}$,
$\sigma^{-1}\left(\frac{1+X}{1+X+X^2}\right)=\frac{(1+X)^2}
{1+X^3+X^4}$, 
$\sigma^{-1}\left(\frac{1+X+X^2}{1+X}\right)=\frac{1+X}
{1+X+X^2}$, $\sigma^{-1}\left(\frac{1+X+X^2}{(1+X)^2}\right)=
\frac{1+X+X^3}{(1+X)^2}$, 
$\sigma^{-1}\left(\frac{1+X+X^2}{(1+X)^3}\right)=
\frac{1+X^3+X^7}{(1+X+X^2)^4}$,
$\sigma^{-1}\left(\frac{1+X+X^2}{(1+X)^4}\right)=
\frac{(1+X+X^3)(1+X^3+X^4)}{(1+X+X^2)^4}$, 

\noindent
$\sigma^{-1}\left(\frac{1+X+X^2}{(1+X)^5}\right)=
\frac{1+X+X^2+X^3+X^4+X^5+X^6+X^{12}+X^{13}}{(1+X+X^2)^7}$, 
$\sigma^{-1}\left(\frac{(1+X+X^2)^2}{1+X}\right)=
\frac{1+X+X^2+X^3+X^4}{(1+X+X^2)^4}$, 
$\sigma^{-1}\left(\frac{(1+X+X^2)^2}{(1+X)^3}\right)=
\frac{(1+X+X^2)(1+X^2+X^5)}{(1+X)^4}$.

\subsection{A few iterations of $\sigma$ and $\sigma^{-1}$ on
  rational fractions in $1+X\mathbb F_2[X]$}\label{sectiterexples}

\subsubsection{Example}
Iterating $\sigma^{-1}$ on $1+X$ yields the following rational fractions
given by their simplest expression, corresponding 
not necessarily to the complete 
factorisation into irreducible polynomials of their numerators and 
enumerators (such a factorisation makes sense when working 
in the multiplicative algebra $\mathbb F_2[[X]]$
and is probably irrelevant for the map $\sigma$, related to the 
shuffle algebra structure
$(\mathbb F_2[[X]],\sh)$).
$$\begin{array}{l}
\sigma^{-1}(1+X)=\frac{1}{1+X}\\
\sigma^{-2}(1+X)=\frac{1}{1+X+X^2}\\
\sigma^{-3}(1+X)=\frac{(1+X)^3}{1+X^3+X^4}\\
\sigma^{-4}(1+X)=\frac{1+X+X^4+X^5+X^7}{1+X^4+X^6+X^7+X^8}\\
\sigma^{-5}(1+X)=\frac{1+X+X^2+X^3+X^4+X^5+X^7+X^8+X^{14}}{1+X^{15}
+X^{16}}\\
\sigma^{-6}(1+X)=\frac{(1+X)^2(1+X+X^2+X^{14}+X^{17}+X^{20}+X^{21}+X^{24}+X^{25}+X^{26}+X^{29})}{1+X^{16}+X^{30}+X^{31}+X^{32}}\\
\end{array}$$

\subsubsection{Example}
Iterating $\sigma^{-1}$ or $\sigma$ on
  $\frac{1+X+X^2}{(1+X)^2}=
1+X+X^3+X^5+X^7+\dots$
yields the following (not necessarily completely factored) results:
$$\begin{array}{l}
\sigma^{-4}\left(\frac{1+X+X^2}{(1+X)^2}\right)=\frac{1+X+X^3+X^5+X^6+X^8+X^9+X^{10}+X^{13}+X^{14}}
{1+X^8+X^{12}+X^{14}+X^{16}}\\
\sigma^{-3}\left(\frac{1+X+X^2}{(1+X)^2}\right)=\frac{1+X+X^2+X^3+X^5}{(1+X^3+X^4)^2}\\
\sigma^{-2}\left(\frac{1+X+X^2}{(1+X)^2}\right)=\frac{(1+X)^3}{(1+X+X^2)^2}\\
\sigma^{-1}\left(\frac{1+X+X^2}{(1+X)^2}\right)=\frac{1+X+X^3}{(1+X)^2}\\
\sigma^{1}\left(\frac{1+X+X^2}{(1+X)^2}\right)=\frac{1+X+X^4}{(1+X)^2}\\
\sigma^{2}\left(\frac{1+X+X^2}{(1+X)^2}\right)=\frac{1+X+X^8}{(1+X)^4}\\
\sigma^{3}\left(\frac{1+X+X^2}{(1+X)^2}\right)=\frac{1+X+X^2+X^4+X^{10}+X^{12}+X^{16}}{(1+X)^8}\\
\sigma^{4}\left(\frac{1+X+X^2}{(1+X)^2}\right)=\frac{1+X+X^3+X^5+X^6+X^{10}+X^{11}+X^{12}+X^{13}+X^{22}+X^{26}+X^{28}+X^{32}}{(1+X)^{16}}\\
\end{array}$$

\subsubsection{Example}
Iterating $\sigma^{-1}$ on
  $\frac{1}{1+X+X^3}$ yields the following 
(not necessarily completely factored) rational fractions:
$$\begin{array}{l}
\sigma^{-3}\left(\frac{1}{1+X+X^3}\right)=\frac{(1+X+X^2+X^4+X^6+X^{12}+X^{15})(1+X^2+X^5+X^6+X^{10}+X^{12}+X^{15})}{1+X^{24}+X^{28}+X^{31}+X^{32}}\\
\sigma^{-2}\left(\frac{1}{1+X+X^3}\right)=\frac{1+X+X^2+X^3+X^5
+X^8+X^{10}+X^{11}+X^{15}}{1+X^8+X^{14}+X^{15}+X^{16}}\\
\sigma^{-1}\left(\frac{1}{1+X+X^3}\right)=\frac{(1+X)^5}{1+X^7+X^8}\\
\sigma^{1}\left(\frac{1}{1+X+X^3}\right)=\frac{1+X+X^2+X^3+X^4}{1+X^2+X^3}\\
\sigma^{2}\left(\frac{1}{1+X+X^3}\right)=\frac{1+X+X^2+X^3+X^4+X^6+X^8}{(1+X^2+X^3)^2}\\
\sigma^{3}\left(\frac{1}{1+X+X^3}\right)=\frac{1+X+X^4+X^5+X^6+X^8+X^9+X^{10}+X^{12}+X^{13}+X^{14}+X^{16}}{(1+X^2+X^3)^4}\\
\sigma^{4}\left(\frac{1}{1+X+X^3}\right)=\frac{P_4}{(1+X^2+X^3)^{8}}\\
\end{array}$$

\begin{rem} Define the degree of a rational fraction
$A\in\mathbb
F_2[[x]]$ as 
$\mathop{deg}(A)=\hbox{max}({\mathop{deg}(f),\mathop{deg}(g)})$ 
if $A=\frac{f}{g}$ with
$f,g\in\mathbb F_2[x]$ without common factor.
For rational $A\in 1+X\mathbb F_2[[X]]$ we have 
$\hbox{lim}_{n\rightarrow\pm \infty}\frac{1}{\vert n\vert}\mathop{log}
(\mathop{deg}(\sigma^n A))=0$ if the orbit of $A$ under $\sigma$ is 
finite. The three examples of Section \ref{sectiterexples} 
suggest that this limit exists (and equals 
$\mathop{log}(2)$) for these examples). It would be
interesting to prove the existence of this limit (or to exhibit a
counterxample) for an arbitrary rational fraction
$A\in 1+{\mathbb F_2}[[X]]$. Since we have
clearly $\mathop{lim}_{n\rightarrow\infty}
\frac{1}{n}\mathop{log}
(\mathop{deg}(\sigma^n(A)))=\mathop{log}(2)$ for $A\in\mathbb F_2[X]$
a polynomial with infinite orbit, one can also ask for 
the existence of values other than $0,\mathop{log}(2)$ for this limit 
which defines obviously an invariant of orbits under the 
bijection $\sigma$ on rational fractions in 
$1+{\mathbb F_2}[[X]]$.
\end{rem}

\subsection{Examples with algebraic series in $1+X\mathbb F_2[[X]]$}
\label{sectalgcplproof}

An algebraic power series $A=\sum_{n=0}^\infty\alpha_nX^n
\in\overline{\mathbb F_2}[[X]]$ 
can be conveniently described by 
a basis of the finite-dimensional vector space 
$\mathcal K(A)$ introduced in Section \ref{sectidalgc}.
More precisely, given a word $\epsilon_1\dots \epsilon_l\in
\{0,1\}^l$ of finite length $l\in\mathbb N$,
we consider the power series 
$$A_{\epsilon_1\dots,\epsilon_l}=\sum_{n=0}^\infty\alpha_{n2^l+
\sum_{j=1}^l\epsilon_j2^{j-1}}X^n.$$
Properties of the Frobenius map imply the identity
$$A_{\epsilon_1\dots,\epsilon_l}=
A_{\epsilon_1\dots,\epsilon_l0}^2+X
A_{\epsilon_1\dots,\epsilon_l1}^2\ .$$
The expression of these identities in terms of a basis for $\mathcal 
K(A)$ gives a fairly compact descriptions for
algebraic series in $\mathbb F_2[[X]]$ as illustrated 
by a few examples below. 

A minimal polynomial
of an algebraic series $A\in\mathbb F_2[[X]]$ can be of degree
$2^{\mathop{dim}(\mathcal K(A))}$ in the variable $A$. 
One can recover such a minimal
polynomial for $A$ by applying an algorithm for Gr\"obner bases to 
the identities described above associated to polynomial relations
in $\mathcal K(A)$ (in terms of a basis or of a generating set).

\subsubsection{Example} 
The preimage $z=\sigma^{-1}(1+\sum_{n=0}^\infty X^{2^n})$ satisfies
the polynomial equation $1+(1+X)z^3=0$.

\subsubsection{Example} Consider the algebraic series
$y=1+\sum_{n=0}^\infty X^{3\cdot 4^n}$ satisfying $y+y^4+X^3=0$
already considered in Remark \ref{remorbinf}. The series
$z=\sigma^{-1}(y)$ satisfies the algebraic equation $1+(1+X^3)z^3$.

\subsubsection{Example} 
Consider the algebraic power series $y=\sum_{n=0}^\infty
X^{2^n-1}=1+X+X^3+X^7+X^{15}+X^{31}+\dots \in\mathbb F_2[[X]]$
satisfying the polynomial equation $1+y+Xy^2=0$.
The formal power series $z=\sigma^{-1}(y)=1+X+X^2+\dots$ 
satisfies the algebraic equation
$$1+X^2+X^3+(1+X)^4z+X(1+X)^4z^2=0$$
and is given by
$$z=\frac{1}{1+X}+X^3\left(\sum_{n=0}^\infty
(\mathop{tm}(n)+\mathop{tm}(n+1))X^{n}\right)^4\in\mathbb F_2[[X]]$$
where $\mathop{tm}\left(\sum_{j=0}\epsilon_j 2^j\right)=
\sum_{j=0}\epsilon_j$ is the Thue-Morse sequence
(see also \cite{AABBJPS} for the sequence
 $n\longmapsto \mathop{tm}(n)+\mathop{tm}(n+1)\pmod 2$).

\begin{rem} For all $n\in\mathbb N$, one can show that  
$\sigma^n(y)=y+P_n(X)$ with $P_n(X)\in\mathbb F_2[X]$ a polynomial where 
$y=\sum_{n=0}^\infty
X^{2^n-1}$. (The series $\sigma^n(y)$ is of course algebraic for 
all $n\in\mathbb N$, see Proposition \ref{propsigmaratalg}.)
\end{rem}

\subsubsection{Example} 
Consider the algebraic power series $y=\sum_{n=0}^\infty
  \mathop{tm}(n+1)
X^n=1+x+x^3+x^6+\dots\in\mathbb F_2[[X]]$ (satisfying
$(1+(1+x)^2y+x(1+x)^3y^2=0$) related to the Thue-Morse
sequence. The preimage $z=\sigma^{-1}(y)$ yields the 
algebraic system of equations
$$\begin{array}{l}
\displaystyle z=z_0^2+Xz_1^2\\
\displaystyle z_0=z_0^2+Xz_{01}^2\\
\displaystyle z_1=z_{10}^2+Xz_{11}^2\\
\displaystyle z_{01}=z_{01}^2+X(z_0+z_{10})^2\\
\displaystyle
z_{10}=z_{10}^2+X(z_0+z_1+z_{01}+z_{11})^2=z_1+X(z_0+z_1+z_{01})^2\\
\displaystyle z_{11}=(z_1+z_{10}+z_{11})^2+X(z_{01}+z_{10}+z_{11})^2=
z_1+(z_1+z_{11})^2+X(z_{01}+z_{10})^2
\end{array}$$

\subsubsection{Example}
Consider the algebraic series
$y=\sigma^{-1}\left(\sum_{n=0}^\infty (\mathop{tm}(n)+
\mathop{tm}(n+1))X^n\right)\in\mathbb F_2[[X]]$ (satisfying 
$1+(1+X)y+X(1+X)y^2=0$).
The preimage $z=\sigma^{-1}(y)\in\mathbb F_2[[X]]$ 
satisfies the algebraic system of equations:
$$\begin{array}{l}
\displaystyle z=z_0^2+Xz_1^2,\\
\displaystyle z_0=z_{00}^2+Xz_0^2,\\
\displaystyle z_1=Xz_{11}^2,\\ 
\displaystyle z_{00}=z_0^2+X(z_0+z_{00})^2,\\
\displaystyle z_{11}=z_{00}^2+X(z_0+z_1+z_{00})^2
\end{array}$$
which, together with the constant terms
$z(0)=z_0(0)=z_{00}(0)=z_{11}(0)=1,z_1(0)=0$,
determines the series $z,z_0=\frac{1}{1+X+X^2},z_1,z_{00}=
\frac{1+X}{1+X+X^2},z_{11}=z+\frac{X^2(1+X)}{(1+X+X^2)^2}$ uniquely.
Eliminating the series $z_0,z_1,z_{00},z_{11}$ by Gr\"obner-basis
techniques yields the  algebraic equation  
$$1+X^2+X^6+X^{10}+X^{11}+X^{12}+X^{15}+(1+X+X^2)^8z+X^3(1+X+X^2)^8z^4=0\
$$
for $z$.

\subsubsection{Example} The series $y=\sum_{n=0}^\infty {3n\choose
  n}X^n\in \mathbb F_2[[X]]$ satisfies the algebraic equation 
$y=1+Xy^3$ (cf. page 423 of \cite{AS}).
Its preimage $z=\sigma^{-1}(y)$ gives rise to the algebraic system
$$\begin{array}{ll}
\displaystyle z=z_0^2+Xz_1^2,\\
\displaystyle z_0=z^2+Xz_{01}^2,\\
\displaystyle z_1=z_{10}^2+Xz_{11}^2,\\
\displaystyle z_{01}=z_{010}^2+Xz_{011}^2,\\
\displaystyle z_{10}=(z_0+z_{010})^2+Xz_{011}^2&
  \displaystyle =z_{01}+z_0^2,\\
\displaystyle z_{11}=z_{010}^2&
  \displaystyle =(1+X)^2z^4+z_1^2,\\
\displaystyle z_{010}=(z+z_{10})^2+X(z+z_{11})^2&
  \displaystyle =(1+X)z^2+z_1,\\
\displaystyle z_{011}=(z+z_{10})^2+X(z_{01}+z_{10})^2\end{array}$$

\section{Other primes}\label{sectothpr} 

There exists an analogue of the quadratic map 
$\sigma:\overline{\mathbb F_2}[[x]]\longrightarrow
\overline{\mathbb F_2}[[x]]$ for $p$ an arbitrary prime. 
It corresponds to the $p-$homogenous form
(still denoted) $\sigma:\overline{\mathbb F_p}[[X]]
\longrightarrow \overline{\mathbb F_p}[[X]]$
defined by 
$$\sigma(A)\equiv\tilde \alpha_0^p+\sum_{n=1}^\infty \tilde \beta_nX^n\pmod p$$
for $A=\sum_{n=0}^\infty \alpha_nX^n\in\overline{\mathbb F_p}[[X]]$
with $\sum_{n=1}^\infty \tilde \beta_nX^n\in X\overline{\mathbb
    Q}[[X]]$
given by the equality
$$\tilde A^{\sh^p}=\tilde \alpha_0^p+p\left(
\sum_{n=1}^\infty \tilde \beta_nX^n\right)$$
for $\tilde A\in\overline{\mathbb Q}[[X]]$ an integral
lift of $A\equiv\tilde A\pmod p$.

The $p-$homogeneous form $\sigma$ restricts to a
bijection of $1+X\overline{\mathbb F_p}[[X]]$ and shares most 
properties holding for $p=2$. In particular, we have:

\begin{prop} \label{propsigmapratalg}
The formal power series $\sigma(A)$ is rational
  (respectively algebraic) if $A\in\overline{\mathbb F_p}[[X]]$
is rational (respectively algebraic).
\end{prop}

\begin{conj} \label{conjXp}
If $A\in 1+\overline{\mathbb F_p}[[X]]$ is rational
  (respectively algebraic) then its preimage $\sigma^{-1}(A)$
is rational (respectively algebraic).
\end{conj}

\subsection{A few examples for $p=3$}

Values of $\sigma^{-1}(A)\in\mathbb F_3[[X]]$
for a few rational $A\in 1+X\mathbb F_3[[X]]$ are:

$\sigma^{-1}(1+X)=\frac{1}{1-X}$, 
$\sigma^{-1}((1+X)^2)=\frac{1-X-X^2}{(1+X)^3}$,
$\sigma^{-1}(\frac{1}{1+X})=\frac{(1+X)^2}{1-X^2+X^3}$,  
$\sigma^{-1}(\frac{1}{(1+X)^2})=\frac{1+X+X^2-X^4+X^5+X^7+X^8}
{(1-X^2+X^3)^3}$,
$\sigma^{-1}(\frac{1+X}{1-X})=\frac{1-X-X^2}
{1-X^2+X^3}$,
$\sigma^{-1}(\frac{1+X}{1+X^2})=\frac{(1-X)^2}
{(1+X)(1-X-X^2)}$.

\subsubsection{Two algebraic examples for $p=3$}
The algebraic series $\sum_{n=0}^\infty X^{3^n-1}=1+X^2+X^8+X^{26}+\dots$
is fixed by $\sigma$.
 
The preimage $z=\sigma^{-1}(1+\sum_{n=0}^\infty X^{3^n})$
satisfies the polynomial equation $(1+X)^3(1-X)z^{13}-1$.
(The power series $y=1+\sum_{n=0}^\infty X^{3^n}\in\mathbb F_3[[X]]$
satisfies the algebraic equation $y=X+y^3$.)

\subsection{A few rational examples for $p=5$}

We give here values of $\sigma^{-1}(A)\in\mathbb F_5[[X]]$
for a few rational $A\in 1+X\mathbb F_5[[X]]$:

$\sigma^{-1}(1+X)=\frac{1}{1-X}$, 
$\sigma^{-1}((1+X)^2)=\frac{(1-X)(1+2X)(1+X+X^2)}{(1-2X)^5}$,
$\sigma^{-1}((1+X)^3)=\frac{(1-2X)(1+2X^2-X^3)}{(1+2X)^5}$,
$\sigma^{-1}(\frac{1}{1+X})=\frac{(1-2X)(1+X-X^2-2X^3)}{1-X^4+X^5}$,
$\sigma^{-1}(\frac{1}{(1+X)^2})=\frac{1-2X+2X^2+2X^4+2X^5-2X^6-2X^7-X^8+X^9+X^{11}-2X^{13}-2X^{14}-2X^{15}-X^{16}+X^{18}+X^{19}-2X^{21}+X^{24}}{(1-X^4+X^5)^5}$,
$\sigma^{-1}(\frac{1+X}{1-X})=\frac{1+2X+X^2+2X^3-2X^4}{1-X^4-2X^5}$,
$\sigma^{-1}(\frac{1+X}{1-2X})=\frac{1-2X-2X^3}{1-X^4+2X^5}$,
$\sigma^{-1}(\frac{1+X}{1+2X})=\frac{1-X-2X^2-X^3-2X^4}{1-X^4+X^5}$.

\section{Power series in free non-commuting variables}
\label{sectionpows}

This and the next section recall a 
few basic and well-known facts concerning
(rational) power series in free non-commuting variables, 
see for instance \cite{St2}, \cite{BR} or a similar book
on the subject. Our terminology,
motivated by \cite{B}, differs however sometimes in the next section.

We denote by $\mathcal X^*$ the free monoid on a set $\mathcal X=
\{X_1,\dots,X_k\}$. We write $1$ for the identity element 
and we use a boldface capital 
$\mathbf X$ for a non-commutative monomial
$\mathbf X=X_{i_1}X_{i_2}\cdots X_{i_l}\in \mathcal X^*$. 
We denote by
$$A=\sum_{\mathbf X\in\mathcal X^*}(A,\mathbf X)\mathbf X\in
\mathbb K\dlangle X_1,\cdots ,X_k\drangle$$
a non-commutative formal power series where
$$\mathcal X^*\ni\mathbf X\longmapsto (A,\mathbf X)\in \mathbb K$$
stands for the coefficient function.

A formal power series 
$A\in\mathbb K\dlangle X_1,\dots,X_k\drangle $ is invertible with respect to
the obvious non-commutative product if and only if it has non-zero 
constant coefficient. 
We denote by 
$\mathfrak m\subset \mathbb K\dlangle X_1,\dots,X_k\drangle $
the maximal ideal consisting of formal power series without
constant coefficient and by 
$\mathbb K^*+\mathfrak m$ the {\it unit-group} of the
algebra $\mathbb K\dlangle X_1,\dots,X_k\drangle $ which is thus
the non-commutative multiplicative
group consisting of all (multiplicatively) invertible elements in
$\mathbb K\dlangle X_1,\dots,X_k\drangle $. 
The unit group is
isomorphic to the direct product $\mathbb K^*\times (1+\mathfrak m)$
where $\mathbb K^*$ is the central subgroup consisting of
non-zero constants and where $1+\mathfrak m$ 
denotes the multiplicative subgroup given by the affine subspace
spanned by power series with
constant coefficient $1$. We have $(1-a)^{-1}=1+\sum_{n=1}^\infty
a^n$ for the multiplicative inverse $(1-a)^{-1}$
of an element $1-a\in 1+\mathfrak m$. 

\subsection{The shuffle algebra}\label{sectshufflealg}

The {\it shuffle-product} $\mathbf X\sh \mathbf X'$ of two 
non-commutative monomials
$\mathbf X,\mathbf X'\in \mathcal X^*$ of degrees 
$a=\mathop{deg}(\mathbf X)$ and $b=\mathop{deg}(\mathbf X')$ 
(for the obvious
grading given by $\mathop{deg}(X_1)=\dots=\mathop{deg}(X_k)=1$)
is the sum
of all ${a+b\choose a}$ monomials of degree $a+b$ obtained by 
``shuffling'' in every possible way the linear factors (elements of
$\mathcal X$) involved in $\mathbf X$ 
with the linear factors of $\mathbf X'$. Such a monomial contribution to 
$\mathbf X\sh \mathbf X'$ can be thought of as a monomial 
of degree $a+b$ whose
linear factors are coloured by two colours with $\mathbf X$ corresponding 
to the product of all linear factors of the first colour and
$\mathbf X'$ corresponding to the product of the remaining 
linear factors. The shuffle product $\mathbf X\sh \mathbf X'$ 
can also be recursively defined by $\mathbf X\sh 1=1\sh \mathbf X=
\mathbf X$ and 
$$(\mathbf X X_s)\sh(\mathbf X' X_t)=(\mathbf X\sh(\mathbf X' X_t))
X_s+((\mathbf X X_s)\sh
\mathbf X')X_t$$ where
$X_s,X_t\in\mathcal X=\{X_1,\dots,X_k\}$ are monomials of degree $1$.

Extending 
the shuffle-product in the obvious way to formal power series endows the
vector space $\mathbb K\dlangle X_1,\dots,X_k\drangle $ with an associative and
commutative algebra structure called the
{\it shuffle-algebra} which has close connections with multiple zeta
values, the algebra of quasi-symmetric functions etc, see eg. \cite{Ho}.
In the case of one variable $X=X_1$ we recover the
definition of Section \ref{sectshufflpr}.

The group $\hbox{GL}_k(\mathbb K)$ acts 
on the vector-space $\mathbb K\dlangle X_1,\dots,
X_k\drangle $ by a linear change of variables. 
This action induces an automorphism of the multiplicative
(non-commutative) algebra or of
the (commutative) shuffle algebra underlying $\mathbb K\dlangle X_1,
\dots,X_k\drangle$. 

Identifying all variables $X_j$ of a formal power series $A
\in\mathbb K\dlangle X_1,\dots,X_k\drangle $ with a common variable 
$X$ yields a
homomorphism of algebras (respectively shuffle-algebras)
from $\mathbb K\dlangle X_1,\dots,X_k\drangle $ into the commutative algebra 
(respectively into the shuffle-algebra) $\mathbb K[[X]]$.

The commutative unit group (set of invertible elements for the
shuffle-product) of the shuffle algebra is given by the set
$\mathbb K^*+\mathfrak m$ and is isomorphic to the direct product
$\mathbb K^*\times (1+\mathfrak m)$. The inverse of an element
$1-a\in1+\mathfrak m$ is given by $\sum_{n=0}^\infty a^{\sh^n}=
1+a+a\sh a+a\sh a\sh a+\dots$, cf. Remark \ref{remamoins1inverse}. 

The following result generalises Proposition \ref{propshgrFp}:

\begin{prop} Over a field of positive characteristic $p$,
the subgroup $1+\mathfrak m$ of the shuffle-group is an
$\mathbb F_p-$vector space of infinite dimension.
\end{prop}

{\bf Proof} Contributions to a $p-$fold shuffle product
$A_1\sh A_2\sh\cdots \sh A_p$ are given by 
monomials with linear factors coloured by 
$p$ colours $\{1,\dots,p\}$
keeping track of their ``origin'' with coefficients given
by the product of the corresponding ``monochromatic'' coefficients
in $A_1,\dots,A_p$. A permutation of 
the colours $\{1,\dots,p\}$ (and in particular, 
a cyclic permutation of all 
colours) leaves such a contribution invariant if
$A_1=\dots=A_p$. Forgetting the colours, coefficients of degree $>0$
in $A^{\sh^p}$ are thus zero in characteristic $p$.
\hfill$\Box$

\section{Rational formal power series}

A formal power series $A$ is {\it rational} if it belongs to
the smallest 
subalgebra in $\mathbb K\dlangle X_1,\dots,X_k\drangle $
which contains the free
associative algebra $\mathbb K\langle
X_1,\dots,X_k\rangle$ of non-commutative polynomials
and intersects the multiplicative unit group of
$\mathbb K\dlangle X_1,\dots,X_k\drangle $ in a subgroup.

The (generalised) {\it Hankel matrix} $H=H(A)$ of
$$A=\sum_{\mathbf X\in\mathcal X^*}(A,\mathbf X)\mathbf X\in
\mathbb K\dlangle X_1,\dots,X_k\drangle $$ 
is the infinite matrix with rows and columns 
indexed by the free monoid $\mathcal X^*$ of monomials and entries
$H_{\mathbf X\mathbf X'}=(A,\mathbf X\mathbf X')$. 
In analogy with the terminology of \cite{B}, 
we call the rank $\mathop{rank}(H)\in\mathbb
N\cup\{\infty\}$ the {\it complexity} of $A$. The row-span,
denoted by $\overline A$,
of $H$ is the {\it recursive closure} of $A$. It corresponds to the 
syntaxic ideal of \cite{BR} and its dimension
$\mathop{dim}(\overline A)$ is the complexity of $A$.

\begin{rem} \label{remdimratfrac}
In the case of one variable, the complexity 
$\mathop{dim}(\overline A)$ of a non-zero 
rational fraction $A=\frac{f}{g}$ with
$f\in\mathbb K[X]$ and $g\in1+X\mathbb K[X]$ is given by
$\mathop{dim}(\overline
A)=\mathop{max}(1+\mathop{deg}(f),\mathop{deg}(g))$.
\end{rem}

Rational series coincide with series of finite complexity by a 
Theorem of Sch\"utzenberger (cf. \cite{BR}, Theorem 1 of page 22).

We call a subspace
$\mathcal A\subset \mathbb K\dlangle X_1,\dots,X_k\drangle $ {\it recursively
  closed} if it contains the recursive closure of all its elements. 

Given a monomial $\mathbf T\in\mathcal X^*$, we
denote by $$\rho(\mathbf T):\mathbb K \dlangle X_1,\dots,X_k\drangle 
\longrightarrow \mathbb K \dlangle X_1,\dots,X_k\drangle $$ 
the linear application which associates to 
$A=\sum_{\mathbf X\in\mathcal X^*}(A,\mathbf X)\mathbf X
\in \mathbb K\dlangle X_1,\dots,X_k\drangle $ the formal power series 
$\rho(\mathbf T)A=\sum_{\mathbf X\in\mathcal X^*}(A,\mathbf X\mathbf
T)\mathbf X$.
We have $\rho(\mathbf T)\rho(\mathbf T')=\rho(
\mathbf T\mathbf T')$. It is easy to check that the set $\{\rho(\mathbf
T)A\}_{\mathbf T\in\mathcal X^*}$ spans the recursive closure 
$\overline A$ of a power series $A$.

\begin{thm} \label{thmratclosed}
We have the inclusion
$$\overline{A\sh B}\subset \overline A\sh \overline B$$
for the shuffle product $A\sh B$ of $A,B\in\mathbb K\dlangle X_1,\dots,
X_k\drangle $.
\end{thm}

\begin{cor} \label{corratclosed}
We have
$$\mathop{dim}(\overline{A\sh B})\leq 
\mathop{dim}(\overline A)\ \mathop{dim}(\overline B)$$
for the shuffle product $A\sh B$ of $A,B\in\mathbb K\dlangle X_1,\dots,
X_k\drangle $.

In particular, shuffle products of rational elements in $\mathbb
K\dlangle X_1,\dots,X_k\drangle $ are rational.
\end{cor}

{\bf Proof of Theorem \ref{thmratclosed}} The shuffle product $A\sh
B$ is clearly contained in the vector space 
$$\overline A\sh \overline B=\{Y\sh Z\ \vert Y\in\overline A,Z\in 
\overline B\}\ .$$
For $Y\in\overline A,Z\in 
\overline B$ and $X_s\in\mathcal X=\{X_1,\dots,X_k\}$, 
the recursive definition of the shuffle product given in Section
\ref{sectshufflealg} shows
$$\rho(X_s)(Y\sh Z)=(\rho(X_s)Y)\sh Z+ Y\sh (\rho(X_s)Z)
\in \overline A\sh Z+Y\sh \overline B\subset \overline A\sh
\overline B$$
and the vector space
$\overline A\sh \overline B$ is thus recursively closed.
\hfill$\Box$

\begin{rem} Similar arguments show that the set of 
rational series is also closed under the ordinary product
(and multiplicative inversion of invertible series),
Hadamard product and composition (where
one considers $A\circ(B_1,\dots,B_k)$ with $A\in\mathbb
K\dlangle X_1,\dots,X_k\drangle $ and $B_1,\dots,B_k\in\mathfrak m\subset 
\mathbb K\dlangle X_1,\dots, X_k\drangle $).
\end{rem}

\begin{rem} The shuffle inverse of a rational element in 
$\mathbb K^*+\mathfrak m$ is in general not rational in characteristic
$0$. An exception is given by geometric progressions
$\frac{1}{1-\sum_{j=1}^k \lambda_j X_j}=\sum_{n=0}^\infty 
\left(\sum_{j=1}^k\lambda_j X_j\right)^n$ since we have
$$\frac{1}{1-\sum_{j=1}^k \lambda_j X_j}\sh
\frac{1}{1-\sum_{j=1}^k \mu_j X_j}=
\frac{1}{1-\sum_{j=1}^k (\lambda_j+\mu_j) X_j}\ .$$
(This identity corresponds to the equality $e^{\lambda X}e^{\mu
  X}=e^{(\lambda+\mu)X}$ 
in the case of a unique variable $X=X_1$, see Remark \ref{remisomalg}.)

By Remark \ref{remshfflegrpnoalgelts}, there are no other such 
elements in $1+\mathfrak m$ in the case of a unique variable $X=X_1$. 
I ignore if the maximal rational shuffle subgroup
of $1+\mathfrak m\subset \mathbb K\dlangle X_1,
\dots,X_k\drangle $ (defined as the set
of all rational elements in $1+\mathfrak m$ with rational inverse
for the shuffle product) contains 
other elements if $k\geq 2$ and if $\mathbb K$ is a suitable field of
characteristic $0$.
\end{rem}

\begin{rem} Any
finite set of rational elements in $\mathbb K\dlangle X_1,
\dots,X_k\drangle$ over a field $\mathbb K$ of positive 
characteristic is included in
a unique minimal finite-dimensional recursively closed subspace
of $\mathbb K\dlangle X_1,\dots,X_k\drangle $ which 
intersects the shuffle group
$\mathbb K^*+\mathfrak m$ in a subgroup.  
\end{rem}

\section{The $p-$homogeneous form  $
\sigma:\overline{\mathbb F_p}\dlangle X_1,\dots,X_k\drangle \longrightarrow
\overline{\mathbb F_p}\dlangle X_1,\dots,X_k\drangle $}\label{sectnonlin}


Considering an integral lift 
$\tilde A=\tilde\alpha+\tilde a\in\overline{\mathbb
  Q}\dlangle X_1,\dots,X_k\drangle $ with coefficients in 
algebraic integers 
of $A=\alpha+a\in\alpha+\mathfrak m\subset \overline{\mathbb
  F_p}\dlangle X_1,\dots,X_k\drangle $, we define $\sigma(A)$ by
the reduction of $\tilde\alpha^p+\tilde b$ modulo $p$ where 
$$\tilde A^{\sh^p}=\tilde \alpha^p +p\tilde b\in
\tilde\alpha^p+\mathfrak m\subset\overline{\mathbb Q}\dlangle X_1,\dots,X_k\drangle 
\ .$$

This definition corresponds to the definition of $\sigma$ given 
in Section \ref{sectothpr} in the case of one variable $X=X_1$.

\begin{prop} \label{propphomXk} One has 
$$\mathop{dim}(\overline{\sigma(A)})\leq 1+{\mathop{dim}(\overline
  A)+p-1
\choose p}$$
for $A\in \overline{\mathbb F_p}\dlangle X_1,\dots,X_k\drangle $.

In particular, $\sigma(A)$ is rational for rational $A\in
\overline{\mathbb F_p}\dlangle X_1,\dots,X_k\drangle $.
\end{prop}

{\bf Proof} It is always possible to choose an integral lift $\tilde A
\in\overline{\mathbb Q}\dlangle X_1,\dots,X_k\drangle $ of $A\in\overline{\mathbb
  F_p}\dlangle X_1,\dots,X_k\drangle $
such that $\mathop{dim}(\overline{\tilde A})=
\mathop{dim}(\overline{A})$. The inclusion
$$\overline{\left(\tilde A^{\sh^p}\right)}\subset \left(
\overline{\tilde A}\right)^{\sh^p}$$
implies then easily the result.\hfill$\Box$

It is easy to show that $\sigma$
induces a bijection on the subset $1+\mathfrak m\subset 
\mathbb K\dlangle X_1,\dots,X_k\drangle $ for 
a field $\mathbb K\subset \overline{\mathbb F_p}$. Computations of a few
examples in $\mathbb F_2\dlangle X_1,X_2\drangle $ suggest:

\begin{conj} \label{conjratX_k} 
The formal power series $\sigma^{-1}(A)$ is rational 
for rational $A\in 1+\mathfrak m\subset \overline{\mathbb F_p}\dlangle X_1,
\dots,X_k\drangle $.
\end{conj}


{\bf Acknowledgements} I thank J-P. Allouche, P. Arnoux, 
M. Brion, A. Pantchichkine, T. Rivoal,
J. Sakarovitch, B. Venkov and J-L. Verger-Gaugry for their interest
and helpful remarks.

\noindent Roland BACHER

\noindent INSTITUT FOURIER

\noindent Laboratoire de Math\'ematiques

\noindent UMR 5582 (UJF-CNRS)

\noindent BP 74

\noindent 38402 St Martin d'H\`eres Cedex (France)
\medskip

\noindent e-mail: Roland.Bacher@ujf-grenoble.fr

\end{document}